\documentclass[12pt]{amsart}
\usepackage{amsmath,amsthm,latexsym,amscd,amsbsy,amssymb,amsfonts,fleqn,leqno,
euscript, pb-diagram}

\numberwithin{equation}{section}

\newtheorem{thm}{Theorem}[section]
\newtheorem{pro}[thm]{Proposition}
\newtheorem{lem}[thm]{Lemma}
\newtheorem{cor}[thm]{Corollary}

\begin{document}


\title[Special embeddings of finite-dimensional compacta in Euclidean spaces]
{Special embeddings of finite-dimensional compacta in Euclidean
spaces }

\author{Semeon Bogatyi}
\address{Faculty of Mechanics and Mathematics, Moscow State University, Vorob'evy
gory, Moscow, 119899 Russia} \email{bogatyi@inbox.ru}
\thanks{The first author was supported by Grants NSH 1562.2008.1
and RFFI 09-01-00741-a.}

\author{Vesko  Valov}
\address{Department of Computer Science and Mathematics, Nipissing University,
100 College Drive, P.O. Box 5002, North Bay, ON, P1B 8L7, Canada}
\email{veskov@nipissingu.ca}
\thanks{Research supported in part by NSERC Grant 261914-08}




\begin{abstract}
If $g$ is a map from a space $X$ into $\mathbb R^m$ and $z\not\in
g(X)$, let $P_{2,1,m}(g,z)$ be the set of all lines
$\Pi^1\subset\mathbb R^m$ containing $z$ such that
$|g^{-1}(\Pi^1)|\geq 2$. We prove that for any $n$-dimensional
metric compactum $X$ the functions $g\colon X\to\mathbb R^m$, where
$m\geq 2n+1$, with $\dim P_{2,1,m}(g,z)\leq 0$ for all $z\not\in
g(X)$ form a dense $G_\delta$-subset of the function space
$C(X,\mathbb R^m)$. A parametric version of the above theorem is
also provided.
\end{abstract}

\maketitle

\markboth{}{Embeddings in Euclidean spaces}



\section{Introduction}
In this paper we assume that all topological spaces are metrizable
and all single-valued maps are continuous.

Everywhere below by $M_{m,d}$ we denote the space of all
$d$-dimensional planes $\Pi^d$ (br., $d$-planes) in $\mathbb R^m$.
If $g$ is a map from a space $X$ into $\mathbb R^m$, $q$ is an
integer and $z\not\in g(X)$, let $P_{q,d,m}(g,z)=\{\Pi^d\in
M_{m,d}:|g^{-1}(\Pi^d)|\geq q{~}\mbox{and}{~}z\in\Pi^d\}$. There is
a metric topology on $M_{m,d}$, see \cite{dnf}, and we consider
$P_{q,d,m}(g,z)$ as a subspace of $M_{m,d}$ with this topology.

One of the results from authors' paper \cite{bv1} states that if $X$
is a metric compactum of dimension $n$ and $m\geq 2n+1$, then the
function space $C(X,\mathbb R^m)$ contains a dense $G_\delta$-subset
of maps $g$ such that the set $\{\Pi^1\in
M_{m,1}:|g^{-1}(\Pi^d)|\geq 2\}$ is at most $2n$-dimensional. The
next theorem provides more information concerning the above result:

\begin{thm}
Let $X$ be a metric compactum of dimension $\leq n$ and $m\geq
2n+1$. Then the maps $g\colon X\to\mathbb R^m$ such that $\dim
P_{2,1,m}(g,z)\leq 0$ for all $z\not\in g(X)$ form a dense
$G_\delta$-subset of $C(X,\mathbb R^m)$.
\end{thm}

Theorem 1.1 admits a parametric version.

\begin{thm}
Let $f\colon X\to Y$ be a perfect $n$-dimensional map between
metrizable spaces with $\dim Y=0$, and $m\geq 2n+1$. Then the maps
$g\colon X\to\mathbb R^m$ such that $\dim
P_{2,1,m}(g|f^{-1}(y),z)\leq 0$ for all restrictions $g|f^{-1}(y)$,
$y\in Y$, and all $z\not\in g(f^{-1}(y)$ form a dense
$G_\delta$-subset of $C(X,\mathbb R^m)$ equipped with the source
limitation topology.
\end{thm}

For any map $g\in C(X,\mathbb R^m)$ and $z\not\in g(X)$ we also
consider the set $D_{2,1,m}(g,z)$ consisting of points
$y=(y_1,y_2)\in(\mathbb R^m)^2$ such that $y_1$ and $y_2$ belong to
a line $\Pi^1\subset\mathbb R^m$ with $z\in\Pi^1$, and there exist
two different points $x_1,x_2\in X$ with $g(x_i)=y_i$, $i=1,2$.

Theorem 1.3 below follows from the proof of Theorem 1.2 by
considering the sets $D_{2,1,m}(g,z)$ instead of $P_{2,1,m}(g,z)$.

\begin{thm}
Let $X, Y$, $f$ and $m$ satisfy the hypotheses of Theorem $1.2$.
Then the maps $g\colon X\to\mathbb R^m$ such that $\dim
D_{2,1,m}(g|f^{-1}(y),z)\leq 0$ for all restrictions $g|f^{-1}(y)$,
$y\in Y$, and all $z\not\in g(f^{-1}(y)$ form a dense
$G_\delta$-subset of $C(X,\mathbb R^m)$.
\end{thm}

Recall that for any metric space $(M,\rho)$ the source limitation
topology on $C(X,M)$ can be describe as follows: the neighborhood
base at a given function $f\in C(X,M)$ consists of the sets
$B_\rho(f,\epsilon)=\{g\in C(X,M):\rho(g,f)<\epsilon\}$, where
$\epsilon:X\to(0,1]$ is any continuous positive functions on $X$.
The symbol $\rho(f,g)<\epsilon$ means that
$\rho(f(x),g(x))<\epsilon(x)$ for all $x\in X$. It is well know that
for metrizable spaces $X$ this topology doesn't depend on the metric
$\rho$ and it has the Baire property provided $M$ is completely
metrizable.


\section{Preliminaries}
We need some preliminary information before proving Theorem 1.1.
Everywhere in this section we suppose that $q, m,d$ are integers
with $0\leq d\leq m$ and $q\geq 1$. Moreover, the Euclidean space
$\mathbb R^m$ is equipped with the standard norm $||.||_m$. We also
suppose that $X$ is a metric compactum and
$\Gamma=\{B_1,B_2,..,B_q\}$ is a disjoint family consisting of $q$
closed subsets of $X$. For any $g\in C(X,\mathbb R^m)$ and $z\not\in
g(X)$ we denote by $P_\Gamma(g,z)$ the set
$$\displaystyle\{\Pi^d\in M_{m,d}:g^{-1}(\Pi^d)\cap
B_i\neq\varnothing {~}\mbox{for
each}{~}i=1,..,q{~}\mbox{and}{~}z\in\Pi^d\}.$$ Now, consider the
open subset $\mathcal R^m_X$ of $C(X,\mathbb R^m)\times\mathbb R^m$
consisting of all pairs $(g,z)$ with $z\not\in g(X)$. Define the
set-valued map
\begin{center}
$\Phi_{\Gamma}\colon \mathcal R^m_X\to M_{m,d},$
$\Phi_{\Gamma}(g,z)=P_\Gamma(g,z)$.
\end{center}

\begin{pro}
$\Phi_{\Gamma}$ is an upper semi-continuous and closed-valued map.
\end{pro}

\begin{proof}
Suppose $(g_0,z_0)\in\mathcal R^m_X$. We need to show that for any
open $W\subset M_{m,d}$ containing $\Phi_{\Gamma}(g_0,z_0)$ there
are neighborhoods $O(g_0)\subset C(X,\mathbb R^m)$ and
$O(z_0)\subset\mathbb R^m$ with $O(g_0)\times O(z_0)\subset\mathcal
R^m_X$ and $\Phi_{\Gamma}(g,z)\subset W$ for all $(g,z)\in
O(g_0)\times O(x_0)$. Assume this is not true. Then there exists a
sequence $\{(g_k,z_k)\}_{k\geq 1}\in\mathcal R^m_X$ converging to
$(g_0,z_0)$ and $\Pi^d_k\in P_\Gamma(g_k,z_k)$ with $\Pi^d_k\not\in
W$, $k\geq 1$. For any $i\leq q$ and $k\geq 1$ there exists a point
$x_k^i\in B_i\cap g_k^{-1}(\Pi^d_k)$. Since $A=\bigcup_{i\leq
q}g_0(B_i)\subset\mathbb R^m$ is compact, we take a closed ball $K$
in $\mathbb R^m$ with center the origin containing $A$ in its
interior. Because every $\Pi^d\in P_\Gamma(g_0,z_0)$ intersects $A$,
we can identify $P_\Gamma(g_0,z_0)$ with $\{\Pi^d\cap K:\Pi^d\in
P_\Gamma(g_0,z_0)\}$ considered as a subspace of $\mathrm{exp}(K)$
(here $\mathrm{exp}(K)$ is the hyperspace of all compact subset of
$K$ equipped with the Vietoris topology).

Because $\{g_k\}_{k\geq 1}$ converges in $C(X,\mathbb R^m)$ to
$g_0$, we can assume that $K$ contains each set $\bigcup_{i\leq
q}g_k(B_i)$, $k\geq 1$. Hence, $g_k(x_k^i)\in K\cap\Pi^d_k$ for all
$i\leq q$ and $k\geq 1$. Therefore, passing to subsequences, we may
suppose that there exist points $x_0^i\in B_i$ $i\leq q$, and a
plane $\Pi^d_0\in M_{m,d}$ such that each sequence $\{x_k^i\}_{k\geq
1}$, $i=1,2,..,q$, converges to $x_0^i$ and $\{\Pi^d_k\cap
K\}_{k\geq 1}$ converges to $\Pi^d_0\cap K$. So,
$\lim\{g_0(x_k^i)\}_{k\geq 1}=g_0(x_0^i)$, $i=1,2,..,q$. Then each
$\{g_k(x_k^i)\}_{k\geq 1}$ also converges to $g_0(x_0^i)$.
Consequently, $g_0(x_0^i)\in\Pi^d_0$ for all $i$. Moreover, since
$z_k\in\Pi^d_k$ for all $k$, we also have $z_0\in\Pi^d_0$. Hence,
$\Pi^d_0\in P_\Gamma(g_0,z_0)$, i.e., $\Pi^d_0\in W$. On the other
hand, $W$ is open in $M_{m,d}$ and $\lim\{\Pi^d_k\cap K\}_{k\geq
1}=\Pi^d_0\cap K$ implies that $\{\Pi^d_k\}_{k\geq 1}$ converges to
$\Pi^d_0$ in $M_{m,d}$. This yields $\Pi^d_k\in W$ for almost all
$k$, a contradiction.

The above arguments also show that $P_\Gamma(g,z)$ is closed in
$M_{m,d}$ for all $(g,z)\in\mathcal R^m_X$. So, $\Phi_{\Gamma,m,d}$
is a closed-valued map.
\end{proof}

Let $X$ and the integers $q,d,m$ be as above. We choose a countable
family $\mathcal B$ of closed subsets of $X$ such that the interiors
of the elements of $\mathcal B$ form a base for the topology of $X$.
Let also $$\mathcal R^m_X(k)=\{(g,z)\in C(X,\mathbb
R^m)\times\mathbb R^m: ||z||_m\leq
k{~}\mbox{and}{~}\rho_m(z,g(X))\geq 1/k,$$ where $\rho_m$ is the
standard Euclidean metric on $\mathbb R^m$ and $k$ an integer. If
$\Gamma\subset\mathcal B$ is a disjoint family of $q$ elements, for
any integers $k,s$ and $\epsilon>0$ we consider the set $\mathcal
H_\Gamma(k,s,\epsilon)$ of all maps $g\in C(X,\mathbb R^m)$ such
that each $P_\Gamma(g,z)$, where $(g,z)\in\mathcal R^m_X(k)$, can be
covered by an open in $M_{m,d}$ family $\omega(g,z)$ satisfying the
following conditions:
\begin{itemize}
\item[(1)] $\mathrm{mesh}(\omega(g,z))<\epsilon$;
\item[(2)] the order of $\omega(g,z)$ is $\leq s$ (i.e., each point from $M_{m,d}$
is contained in at most $s+1$ elements of $\omega(g,z)$).
\end{itemize}

\begin{pro}
Any $\mathcal H_\Gamma(k,s,\epsilon)$ is open in $C(X,\mathbb R^m)$.
\end{pro}

\begin{proof}
Assume $g_0\in\mathcal H_\Gamma (k,s,\epsilon)$. For any
$(g_0,z)\in\mathcal R^m_X(k)$ let
$W(g_0,z)=\bigcup\{U:U\in\omega(g_0,z)\}$. Obviously, we have
$(g_0,z)\in\mathcal R^m_X(k)$ if and only if $z$ belongs to the
compact set $B(g_0)=\{z\in\mathcal R^m:||z||_m\leq
k{~}\mbox{and}{~}\rho_m(z,g_0(X))\geq 1/k\}$. Hence,
$P_\Gamma(g_0,z)\subset W(g_0,z)$ for every $z\in B(g_0)$. According
to Proposition 2.1, for any such $z$ there exists an open
neighborhood $O(z)\subset\mathbb R^m\backslash g_0(X)$ such that
$P_\Gamma(g_0,u)\subset W(g_0,z)$ for all $u\in O(z)$. Next, shrink
each $O(z)$ to an open set $V(z)$ such that $z\in
V(z)\subset\overline{V(z)}\subset O(z)$. Then $\{V(z):z\in B(g_0)\}$
is an open cover of $B(g_0)$ and we choose a finite subcover
$\{V(z_j):j=1,2,.,p\}$. Let $\eta$ be the distance between $B(g_0)$
and $\mathbb R^m\backslash V$, where $V=\bigcup_{j=1}^{j=p}V(z_j)$,
and $A(z)=\{j:z\in O(z_j)\}$, $z\in O=\bigcup_{j=1}^{j=p}O(z_j)$.
Choosing smaller neighborhoods $V(z_j)$, if necessarily, we may
assume that $\eta<1/k$. According to the choice of $O(z_j)$, we have
$$P_\Gamma(g_0,z)\subset W(g_0,z_j){~}\mbox{for any}{~}z\in O{~}\mbox{and}{~}j\in A(z)\leqno{(3)}.$$

{\em Claim $1$. Let $g\in O(g_0,\eta)$ and $\rho_m(z,g(X))\geq 1/k$,
where $O(g_0,\eta)$ consists of all $g\in C(X,\mathbb R^m)$ such
that $\rho_m(g_0(x),g(x))<\eta$ for all $x\in X$. Then
$\rho_m(z,g_0(X))\geq (1/k)-\eta$ and $z\in\overline{V}\subset O$}.
\smallskip

Indeed, both $\rho_m(z,g_0(X))<(1/k)-\eta$ and $g\in O(g_0,\eta)$
imply the existence of $x\in X$ with $\rho_m(g_0(x),g(x))<1/k$ which
contradicts $\rho_m(z,g(X))\geq 1/k$. So, for every $z$ satisfying
the hypotheses of Claim 1, we have $\rho_m(z,g_0(X))\geq
(1/k)-\eta$. This yields $z\in\overline{V}\subset O$.

Each $W(g_0,z_j)$ is the union of an open family in $M_{m,d}$ of
order $\leq s$ and $\mathrm{mesh}<\epsilon$. Thus, according to
Claim 1, it suffices to show the next claim.

\smallskip
{\em Claim $2$. There exists a neighborhood $O(g_0)\subset
O(g_0,\eta)$ of $g_0$ satisfying the following condition: for any
$z\in\overline{V}$ with $\rho_m(z,g_0(X))\geq (1/k)-\eta$ there
exists $j\in A(z)$ such that $P_\Gamma(g,z)\subset W(g_0,z_j)$
whenever $g\in O(g_0)$.}
\smallskip

 Suppose the conclusion of Claim 2 is not true. Then for every
$p\geq 1$ there exists a map $g_p\in O(g_0,\eta)$ with
$\rho_m(g_0(x),g_p(x))<1/p$ for all $x\in X$, a point
$z_p\in\overline{V}$ with $\rho_m(z_p,g_0(X))\geq (1/k)-\eta$, and
planes $$\Pi^d_p\in
P_\Gamma(g_p,z_p)\backslash\bigcup\{W(g_0,z_j):j\in
A(z_p)\}.\leqno{(4)}$$ Passing to subsequences, we may assume that
the sequence $\{z_p\}_{p\geq 1}$ converges to a point
$z_0\in\overline{V}$ and $\{\Pi^d_p\}_{p\geq 1}$ converges in
$M_{m,d}$ to a $d$-plane $\Pi^d_0$. Obviously,
$\rho_m(z_0,g_0(X))\geq (1/k)-\eta$. Since $z_p\in\Pi^d_p$, we also
have $z_0\in\Pi^d_0$. As in the proof of Proposition 2.1, we can see
that $g_0^{-1}(\Pi^d_0)$ meets each element of $\Gamma$.
Consequently, $\Pi^d_0\in P_\Gamma(g_0,z_0)$. So, by (3),
$\Pi^d_0\in\bigcap\{W(g_0,z_j):j\in A(z_0)\}$. This implies that
$\Pi^d_p\in\bigcap\{W(g_0,z_j):j\in A(z_0)\}$ for almost all $p$. On
the other hand, since $\lim z_p=z_0$, there exists $p_0$ such that
$A(z_0)\subset A(z_p)$ for all $p\geq p_0$. So, by (4),
$\Pi^d_p\not\in\bigcup\{W(g_0,z_j):j\in A(z_0)\}$ when $p\geq p_0$,
a contradiction.
\end{proof}

\begin{cor}
All maps $g\in C(X,\mathbb R^m)$ such that $\dim P_{q,d,m}(g,z)\leq
s$ for all $z\not\in g(X)$ form a $G_\delta$-subset $\mathcal
H_X(q,d,m,s)$ of $C(X,\mathbb R^m)$.
\end{cor}

\begin{proof}
It easily seen that each $g\in C(X,\mathbb R^m)$ and $z\not\in g(X)$
we have $P_{q,d,m}(g,z)=\bigcup\{P_\Gamma(g,z):\Gamma\subset\mathcal
B {~}\mbox{is disjoint and}{~}|\Gamma|=q\}$. Moreover, since
$P_\Gamma(g,z)$ are closed in $M_{m,d}$ (by Proposition 2.1), we
have $\dim P_{q,d,m}(g,z)\leq s$ if and only if $\dim
P_\Gamma(g,z)\leq s$ for all $\Gamma$. This implies that $\mathcal
H_X(q,d,m,s)$ is the intersection of the sets $\mathcal
H_\Gamma(k,s,1/p)$, where $k,p\geq 1$ are integers and
$\Gamma\subset\mathcal B$ is a disjoint family of $q$ elements.
\end{proof}

\section{Proof of Theorems 1.1 and 1.2}
Recall that a real number $v$ is called algebraically dependent on
the real numbers $u_1,..,u_k$ if $v$ satisfies the equation
$p_0(u)+p_1(u)v+...+p_n(u)v^n=0$, where $p_0(u),..,p_n(u)$ are
polynomials in $u_1,..,u_k$ with rational coefficients, not all of
them 0. A finite set of real numbers is {\em algebraically
independent} if none of them depends algebraically on the others.
The idea to use algebraically independent sets for proving general
position theorems was originated by Roberts in \cite{r}. This idea
was also applied by Berkowitz and Roy in \cite{br}. A proof of the
Berkowitz-Roy main theorem from \cite{br} was provided by Goodsell
in \cite[Theorem A.1]{g2} (see \cite[Corollary 1.2]{bv} for a
generalization of the Berkowitz-Roy theorem and \cite{g1} for
another application of this theorem). Let us note that any finitely
many points in an Euclidean space $\mathbb R^n$ whose set of
coordinates is algebraically independent are in general position.

\textit{Proof of Theorem $1.1$.} We have to show that the set
$\mathcal H_X(2,1,m,0)$ of all maps  $g\in C(X,\mathbb R^m)$ such
that $\dim P_{2,1,m}(g,z)\leq 0$ for all $z\not\in g(X)$ is dense
and $G_\delta$ in $C(X,\mathbb R^m)$. According to Corollary 2.3,
this set is $G_\delta$. So, it remains to show it is also dense in
$C(X,\mathbb R^m)$. Fix a countable family $\mathcal B$ of closed
subsets of $X$ such that the interiors of its elements is a base for
$X$. Since $\mathcal H_X(2,1,m,0)$ is the intersection of the open
family $$\{\mathcal H_\Gamma (k,0,1/p):\Gamma\subset\mathcal B
{~}\mbox{is disjoint with}{~}|\Gamma|=2{~}\mbox{and}{~}k,p\geq 1\}$$
(see the proof of Corollary 2.3), it suffices to show that each
$\mathcal H_\Gamma (k,0,\epsilon)$ is dense in $C(X,\mathbb R^m)$.
Recall that $\mathcal H_\Gamma (k,0,\epsilon)$ consists of all maps
$g\in C(X,\mathbb R^m)$ such that $P_\Gamma(g,z)$ can be covered by
a disjoint open in $M_{m,1}$ family $\omega$ with
$\mathrm{mesh}(\omega)<\epsilon$ for every map $g$ and every point
$z\in\mathbb R^m$ satisfying the following conditions: $||z||_m\leq
k$ and $\rho_m(z,g(X))\geq 1/k$.

To prove that each $\mathcal H_\Gamma (k,0,\epsilon)$ is dense in
$C(X,\mathbb R^m)$, observe that any map $g\in C(X,\mathbb R^m)$ can
be approximated by maps $f=h\circ l$ with $l\colon X\to K$ and
$h\colon K\to\mathbb R^m$, where $K$ is a finite polyhedron of
dimension $\leq n$. Actually, $K$ can be supposed to be a nerve of a
finite open cover $\beta$ of $X$. Moreover, if we choose $\beta$
such that any its element meets at most one element of
$\Gamma=\{B_1,B_2\}$, then we have $l(B_1)\cap l(B_2)=\varnothing$.
Further, taking sufficiently small barycentric subdivision of $K$,
we can find disjoint subpolyhedra $K_i$ of $K$ with $l(B_i)\subset
K_i$, $i=1,2$. Obviously, for any $z\not\in h(l(X))$ the set
$P_\Gamma(h\circ l,z)$ is contained in $P_\Lambda(h,z)=\{\Pi^1\in
M_{m,1}:h^{-1}(\Pi^1)\cap K_i\neq\varnothing,
i=1,2{~}\mbox{and}{~}z\in\Pi^1\}$, where $\Lambda=\{K_1,K_2\}$.
Therefore, the density of $\mathcal H_\Gamma (k,0,\epsilon)$ in
$C(X,\mathbb R^m)$ is reduced to show that the maps $h\in
C(K,\mathbb R^m)$ such that any $P_\Lambda(h,z)$, $z\not\in h(K)$,
admits a disjoint open cover in $M_{m,1}$ of
$\mathrm{mesh}<\epsilon$ form a dense subset of $C(K,\mathbb R^m)$.
And this follows the next proposition.

\begin{pro}
Let $K_i$, $i=1,2$, be disjoint $n$-dimensional subpolyhedra of a
finite polyhedron $K$. Then the maps $h\in C(K,\mathbb R^m)$ such
that for any $z\not\in h(K)$ the set $\{\Pi^1\in M_{m,1}:
h^{-1}(\Pi^1)\cap K_i\neq\varnothing,
i=1,2{~}\mbox{and}{~}z\in\Pi^1\}$ is of dimension $\leq 0$
 form a dense subset of $C(K,\mathbb R^m)$.
\end{pro}

\begin{proof}
Let $h_0\in C(K,\mathbb R^m)$ and $\delta>0$. We take a subdivision
of $K$ such that $\mathrm{diam}h_0(\sigma)<\delta/2$ for all
simplexes $\sigma$. Let $K^{(0)}=\{a_1,a_2,...,a_t\}$ be the
vertexes of $K$ and $v_j=h_0(a_j)$, $j=1,..,t$. Then, by \cite{br},
there are points $b_j\in\mathbb R^m$ such that the distance between
$v_j$ and $b_j$ is $<\delta/2$ for each $j$ and the coordinates of
all $b_j$, $j=1,..,t$, form an algebraically independent set. Define
a map $h\colon K\to\mathbb R^m$ by $h(a_j)=b_j$ and $h$ is linear on
every simplex of $K$. It is easily seen that $h$ is $\delta$-close
to $h_0$. Without loss of generality, we may suppose that $K_1$ and
$K_2$ are two $n$-dimensional simplexes. Then each $h(K_i)$ is also
an $n$-dimensional simplex in $\mathbb R^m$ generating a plane
$\Pi^n_i\in M_{m,n}$. Since the coordinates of the points
$\{b_j:j=1,..,t\}$ form an algebraically independent set, the planes
$\Pi^n_1$ and $\Pi^n_2$ are skew. Suppose $z\not\in h(K)$. If
$z\in\Pi^n_1$ or $z\in\Pi^n_2$, then there is no line
$\Pi^1\subset\mathbb R^m$ which contains $z$ and meets both $h(K_1)$
and $h(K_2)$. Suppose $z\not\in\Pi^n_1\cup\Pi^n_2$. According to
\cite[Corollary 3.8]{bv1}, there exists at most one line
$\Pi^1\subset\mathbb R^m$ containing $z$ such that $\Pi^1\cap
h(K_i)\neq\varnothing$, $i=1,2$. Hence, for any $z\not\in h(K)$ the
set $\{\Pi^1\in M_{m,1}:h^{-1}(\Pi^1)\cap K_i\neq\varnothing,
i=1,2{~}\mbox{and}{~}z\in\Pi^1\}$ is finite.
\end{proof}

\textit{Proof of Theorem $1.2$.} We fix a metric $d$ generating the
topology of $X$ and for any $g\in C(X,\mathbb R^m)$, $y\in Y$,
$\eta>0$ and $z\not\in g(f^{-1}(y))$ let $\displaystyle
P^\eta(g,y,z)$ be the set of all $\Pi^1\in M_{m,1}$ such that
$z\in\Pi^1$ and there exist two points $x^1,x^2\in g^{-1}(\Pi^1)\cap
f^{-1}(y)$ with $d(x^1,x^2)\geq\eta$ . Obviously,
$$\displaystyle
P_{2,1,m}(g|f^{-1}(y),z)=\bigcup_{k=1}^{\infty}\{P^{1/k}(g,y,z){~}\mbox{for
any}{~} z\not\in g(f^{-1}(y))\}.\leqno{(5)}$$

\smallskip
\textit{Claim $3$. Each $\displaystyle P^\eta(g,y,z)$ is closed in
$\displaystyle P_{2,1,m}(g|f^{-1}(y),z)$}.
\smallskip

The proof of Claim 3 follows the arguments from the proof of
Proposition 2.1.

Now, for $k\geq 1$ and $y\in Y$ consider the set $$\displaystyle
B_g(y,k)=\{z\in\mathbb R^m:||z||_m\leq
k{~}\mbox{and}{~}\rho_m(z,g(f^{-1}(y)))\geq 1/k\}.$$ Next, let
$\displaystyle\mathcal P_\epsilon^\eta(y,k)$ be the set of all maps
$g\in C(X,\mathbb R^m)$ such that for each $z\in B_g(y,k)$ the set
$\displaystyle P^\eta(g,y,z)$ can be covered by a disjoint open in
$M_{m,1}$ family of $\mathrm{mesh}<\epsilon$.  If $F\subset Y$, we
consider the set $\displaystyle\mathcal
P_\epsilon^\eta(F,k)=\bigcap_{y\in F}\mathcal P_\epsilon^\eta(y,k)$.
Obviously the intersection of all $\displaystyle\mathcal
P_{1/s}^\eta(Y,k)$, $s\geq 1$, is the set
$$\displaystyle\mathcal P^\eta(Y,k)=\{g\in C(X,\mathbb R^m):\dim
P^\eta(g,y,z)\leq 0,y\in Y,z\in B_g(y,k)\}.$$ It follows from (5)
that the set $\displaystyle\bigcap_{k,s=1}^{\infty}\mathcal
P^{1/s}(Y,k)$ coincides with the set
$$\displaystyle\mathcal P=\{g\in
C(X,\mathbb R^m):\dim P_{2,1,m}(g|f^{-1}(y),z)\leq 0, y\in Y,
z\not\in g(f^{-1}(y)).$$ So, in order to show that $\mathcal P$ is
dense and $G_\delta$ in $C(X,\mathbb R^m)$, it suffices to show that
each $\displaystyle\mathcal P_\epsilon^\eta(Y,k)$ is open and dense
in $C(X,\mathbb R^m)$.

We are going first to show that any $\displaystyle\mathcal
P_\epsilon^\eta(Y,k)$ is open in $C(X,\mathbb R^m)$. This can be
done following the arguments from \cite[Proposition 5.3]{bv1} using
the next lemma instead of \cite[Lemma 5.2]{bv1}.

\begin{lem}
Let $\displaystyle g_0\in\mathcal P_\epsilon^\eta(y_0,k)$ for some
$y_0\in Y$. Then there exists a neighborhood $V$ of $y_0$ in $Y$ and
$\delta>0$ such that $\displaystyle g\in\mathcal
P_\epsilon^\eta(V,k)$ for all $g\in C(X,\mathbb R^m)$ such that the
restrictions $g|f^{-1}(V)$ and $g_0|f^{-1}(V)$ are $\delta$-close.
\end{lem}

\begin{proof}
Assume the conclusion of Lemma 3.2 doesn't hold and use the
arguments from the proof of Propositions 2.1 and 2.2 to obtain a
contradiction.
\end{proof}

The next proposition completes the proof of Theorem 1.2.

\begin{pro}Any set $\displaystyle\mathcal
P_\epsilon^\eta(Y,k)$ is dense in $C(X,\mathbb R^m)$ with respect to
the source limitation topology.
\end{pro}

\begin{proof}
We modify the arguments from the proof of \cite[Proposition
5.4]{bv1}. Let $g\in C(X,\mathbb R^m)$ and $\delta\in C(X,(0,1])$.
We are going to find $h\in\displaystyle\mathcal
P_\epsilon^\eta(Y,k)$ such that $\rho(g(x),h(x))<\delta(x)$ for all
$x\in X$. By \cite[Proposition 4]{bv2}, $g$ can be supposed to be
simplicially factorizable. This means that there exists a simplicial
complex $D$ and maps $g_D\colon X\to D$, $g^D\colon D\to M$ with
$g=g^D\circ g_D$. Following the proof of \cite[Proposition
3.4]{bv3}, we can find an open cover $\mathcal U$ of $X$, simplicial
complexes $N, L$ and maps $\alpha:X\to N$, $\beta:Y\to L$, $p\colon
N\to L$, $\varphi\colon N\to\mathbb R^m$ and $\delta_1\colon N\to
(0,1]$ satisfying the following conditions, where
$h'=\varphi\circ\alpha$:

\begin{itemize}
\item $\alpha$ is an $\mathcal U$-map and for any $x_1,x_2\in X$ with $d(x_1,x_2)\geq\eta$ we have
$\alpha(x_1)\neq\alpha(x_2)$;
\item $\beta\circ f=p\circ\alpha$;
\item $p$ is a perfect $PL$-map with $\dim p\leq n$ and $\dim L=0$;
\item $h'$ is $(\delta/2)$-close to $g$;
\item $\delta_1\circ\alpha\leq\delta$.
\end{itemize}
 So, we have the following commutative diagram:

 \begin{picture}(120,95)(-100,0)
\put(30,10){$L$} \put(0,30){$Y$} \put(12,28){\vector(3,-2){18}}
\put(14,14){\small $\beta$} \put(1,70){$X$}
\put(5,66){\vector(0,-1){25}} \put(-1,53){\small $f$}
\put(11,73){\vector(1,0){45}} \put(30,77){\small $h'$}
\put(12,68){\vector(3,-2){18}} \put(15,56){\small $\alpha$}
\put(31,50){$N$} \put(35,46){\vector(0,-1){25}}
 \put(37,33){\small $p$}
\put(46,58){\vector(4,3){13}} \put(44,64){\small $\varphi$}
 \put(60,70){$\mathbb R^m$}
\end{picture}

Since $L$ is a $0$-dimensional simplicial complex and $p$ is a
perfect $PL$-map, $N$ is a discrete union of the finite complexes
$K_l=p^{-1}(l)$, $l\in L$. Because $\dim p\leq n$, $\dim K_l\leq n$,
$l\in L$. Applying Theorem 1.1 to each complex $K_l$, we can find a
map $\varphi_1\colon N\to\mathbb R^m$ such that for any $l\in L$ and
$z\not\in\varphi_1(p^{-1}(l))$ we have $\displaystyle\dim
P_{2,1,m}(\varphi_1|p^{-1}(l),z)\leq 0$ and $\varphi_1|p^{-1}(l)$ is
$\theta_l$-close to $\varphi|p^{-1}(l)$, where
$\theta_l=\min\{\delta_1(u):u\in p^{-1}(l)\}$. Moreover, the map
$h=\varphi_1\circ\alpha$ is $\delta$-close to $g$. We claim that
$h\in\displaystyle\mathcal P_\epsilon^\eta(Y,k)$. Indeed, let $y\in
Y$ and $z\in B_h(y,k)$. If $\Pi^1\in\displaystyle P^\eta(h,y,z)$,
then there exist two points $x^i\in h^{-1}(\Pi^1)\cap f^{-1}(y)$,
$i=1,2$, with $d(x^1,x^2)\geq\eta$. According to the choice of the
cover $\mathcal U$, we have $\alpha(x^1)\neq\alpha(x^2)$. Since
$\varphi_1^{-1}(\Pi^1)\cap p^{-1}(\beta(y))$ contains the points
$\alpha(x^i)$, $i=1,2$, we obtain that $\displaystyle\Pi^1\in
P_{2,1,m}(\varphi_1|p^{-1}(\beta(y)),z)$. Thus, we established the
inclusion $\displaystyle P^\eta(h,y,z)\subset
P_{2,1,m}(\varphi_1|p^{-1}(\beta(y)),z)$ which implies
$\displaystyle \dim P^\eta(h,y,z)\leq 0$ for every $y\in Y$ and
$z\in B_h(y,k)$. Consequently, $h\in\displaystyle\mathcal
P_\epsilon^\eta(Y,k)$.
\end{proof}

\textbf{Acknowledgments.} The results from this paper were obtained
during the second author's visit of Department of Computer Science
and Mathematics (COMA), Nipissing University in August 2010. He
acknowledges COMA for the support and hospitality.


\end{document}